\documentclass[11pt, reqno]{amsart}
\usepackage[margin=3.6cm]{geometry}
\usepackage{array,booktabs,tabularx}
\usepackage{tikz,tkz-euclide,pstricks}
\usetikzlibrary{arrows,calc,patterns,backgrounds}
\usepackage{graphicx}
\usepackage{amssymb}
\usepackage{enumerate}
\usepackage{nicefrac}
\usepackage{color}
\usepackage{hyperref}
\usepackage{esint}
\usepackage{mathtools}
\usepackage{dsfont}
\usepackage{ esint }
\usepackage[show]{ed}
\usepackage{comment}

\usepackage{enumitem}
\usepackage{soul}
\usepackage{needspace}
\usepackage[normalem]{ulem}
\definecolor{light-gray}{gray}{0.8}
\definecolor{light-blue}{rgb}{0.53,.8,98}
\definecolor{green1}{RGB}{50,180,50}

\newtheorem{theorem}{Theorem}
\newtheorem{corollary}[theorem]{Corollary}

\theoremstyle{definition}

\newcommand{\R}{{\mathbb R}}

\newcommand{\x}{\ensuremath{\times}}

\newcommand{\mc}[1]{\mathcal{#1}}
\newcommand{\re}{\mathbb{R}}

\renewcommand{\red}[1]{{\color{red}{#1}}}

\newcommand{\SNH}{S\!N^*\!H}

\newcommand{\sig}{\sigma_{_{\!\!S\!N^*\!H}}}

\newcommand{\SM}{S^*\!M}
\newcommand{\TM}{T^*\!M}

\newcommand{\Tj}{\Lambda_{\rho_j}^\tau(R(h))}
\newcommand{\LambdaH}{\Lambda^\tau_{_{\!\Sigma_{H\!,p}}}}

\newcommand{\LM}{{_{\!L^2(M)}}}

\newcommand{\class}{\delta}

\def\XXint#1#2#3{{\setbox0=\hbox{$#1{#2#3}{\int}$} \vcenter{\hbox{$#2#3$}}\kern-.5\wd0}}

\DeclareMathOperator{\supp}{supp}

\newcommand{\e}{\varepsilon}
\newcommand{\Tinj}{\tau_{_{\!\text{inj}H}}}

\numberwithin{equation}{section}
\numberwithin{lemma}{section}

\title[]{A microlocal approach to eigenfunction concentration}

\author{Jeffrey Galkowski}
\address{Department of Mathematics, Northeastern University, Boston, MA, USA}
\email{jeffrey.galkowski@northeastern.edu }

\date{}

\begin{document}

\begin{abstract}
We describe a new approach to understanding averages of high energy Laplace eigenfunctions, $u_h$, over submanifolds,
$$
\Big|\int _H u_hd\sigma_H\Big|
$$
where $H\subset M$ is a submanifold and $\sigma_H$ the induced by the Riemannian metric on $M$. This approach can be applied uniformly to submanifolds of codimension $1\leq k\leq n$ and in particular, gives a new approach to understanding $\|u_h\|_{L^\infty(M)}$. The method, developed in~\cite{GT,Gdefect,CGT,CG17,CG18}, relies on estimating averages by the behavior of $u_h$ microlocally near the conormal bundle to $H$. By doing this, we are able to obtain quantitative improvements on eigenfunction averages under certain uniform non-recurrent conditions on the conormal directions to $H$. In particular, we do not require any global assumptions on the manifold $(M,g)$.

\end{abstract}

\maketitle

\section{Introduction}

In this note, we describe a new approach to understanding concentration properties of high energy eigenfunctions. Although the methods in~\cite{GT,Gdefect,CGT,CG17,CG18} (on which this note is based) apply to the quasimodes of a wide variety of pseudodifferential operators, we focus on the case of the Laplacian on a Riemannian manifold $(M,g)$ of dimension $n$ and consider only eigenfunctions i.e. solutions to
\begin{equation}
\label{e:laplace}
(-h^2\Delta_g-1)u_h=0
\end{equation}
for concreteness. Consider a submanifold $H\subset M$. We are interested in averages of the form 
 $$
 \int_H u_hd\sigma_H
 $$
 where $\sigma_H$ denotes the volume measure induced on $H$ from $M$. We note that, when $H=\{x\}$ is a point in $M$, this average is given precisely by $u_h(x)$. Thus, using our methods we are able to obtain control on $L^\infty$ norms. We do not give the details of many proofs in this note, instead referring to the relevant papers. We review some of the previously existing results, state the new theorems, and describe the ideas central to the proofs. 

Since the middle of the twentieth century~\cite{Ava,Lev,Ho68} many authors have been interested in the growth of eigenfunctions for self-adjoint elliptic operators. In particular, they prove that a solution to~\eqref{e:laplace} satisfies, 
\begin{equation}
\label{e:Linf}
\|u_h\|_{L^\infty(M)}\leq Ch^{\frac{1-n}{2}}\|u_h\|_{L^2(M)}.
\end{equation}

If one considers the case of $(M,g)=(S^2,g_{\text{round}})$, the sphere with the round metric, it is not hard to construct the family of zonal harmonics, $Z_{h}$, with the property that 
$$
ch^{-\frac{1}{2}}\|Z_{h}\|_{L^2(S^2)}\leq \|Z_{h}\|_{L^\infty(S^2)}\leq Ch^{-\frac{1}{2}}\|Z_{h}\|_{L^2(S^2)},\qquad (-h^2\Delta_{S^2}-1)Z_h=0,
$$
and hence that the estimate~\eqref{e:Linf} cannot be improved on a general manifold. Because of this, it is natural to try to understand situations in which~\eqref{e:laplace} is sharp. 
It is also interesting to think of the question of $L^\infty$ norms as averages over points and to generalize that question to averages over submanifolds $H\subset M$. While it is a more recent line of inquiry than that of $L^\infty$ bounds, it dates at least to the early 1980's~\cite{Hej,Good}. The analog of~\eqref{e:laplace} was proved in~\cite{Zel}, where the author shows that if $H$ has codimension $k$, then 
\begin{equation}
\label{e:average}
\Big|\int_H u_hd\sigma_H\Big|\leq Ch^{\frac{1-k}{2}}\|u_h\|_{L^2(M)}.
\end{equation}
Again, for every $1\leq k \leq n$, there are examples on the sphere of dimension $n$ which saturate the estimate~\eqref{e:average} and it is natural to characterize situations in which~\eqref{e:average} can be saturated.

\subsection{A review of previous $L^\infty$ results}
Before we can state the results on a general manifold $M$, we need some concepts from geometry. Let $\TM$ denote the cotangent bundle to $M$, $H\subset M$ a submanifold with conormal bundle $N^*\!H$, and $\SNH$, the unit conormal bundle to $H$,
$$
\SNH:=\big\{(x,\xi)\in N^*\!H\mid |\xi|_{g(x)}=1\big\},
$$
where $|\cdot |_{g}$ denotes the metric induced on $\TM$ by $g$. Note that $S\!N^*\!\{x\}=S^*_xM$ where $S^*_xM:=T^*_xM\cap \SM$. Next, let $G^t:\SM\to \SM$ denote the geodesic flow. 

We define the \emph{first return time} $T_H:\SNH\to [0,\infty]$ by 
$$
T_H(x,\xi):=\inf\{t>0\mid G^t(x,\xi)\in \SNH\}.
$$
We then define the \emph{loop set of $H$}, $\mc{L}_H:=T_H^{-1}([0,\infty))$. Finally, we let $\sig$ be the volume induced on $\SNH$ by the Liouville measure on $\TM$. 
\begin{theorem}[\cite{SoggeZelditch}]
\label{t:loops}
Let $(M,g)$ be a Riemannian manifold of dimension $n$. Suppose that $x\in M$ and $\sigma_{_{\!S^*_xM}}(\mc{L}_x)=0$. Then for $u_h$ solving~\eqref{e:laplace}
$$
|u_h(x)|=o\big(h^{\frac{1-n}{2}}\|u_h\|_{L^2(M)}\big).
$$
\end{theorem}

We define the first return map $\eta_H:\mc{L}_H\to \SNH$ by
$$
\eta_H(x,\xi):=G^{T_H(x,\xi)}(x,\xi).
$$
and let 
$$
\mc{L}_H^{\pm n}:=\bigcap_{k=0}^n\eta_H^{\pm k}(\SNH),\qquad \mc{L}_H^\infty:=\bigcap_n\mc{L}_H^n.
$$
Next, define the \emph{recurrent set of $H$}
$$
\mc{R}_H:=\Big\{(x,\xi) \in \mc{L}_H^\infty \mid (x,\xi) \in \Big[\bigcap_{n>0}\overline{\bigcup_{k\geq n}\eta_H^k(x,\xi)}\,\Big]\bigcap \Big[ \bigcap_{n>0}\overline{\bigcup_{k\geq n}\eta_H^{-k}(x,\xi)}\,\Big]\,\Big\}.
$$
In~\cite{SoggeTothZelditch}, the authors show that $\mc{L}_x$ can be replaced by $\mc{R}_x$ in Theorem~\ref{t:loops}. Finally, in~\cite{SZ16I, SZ16II}, the authors obtain still more restrictive assumptions in the case that $(M,g)$ is real analytic. In fact, in the case of a real analytic surface, they are able to verify the conjecture that one can replace $\mc{R}_x$ by the set of directions $\xi$ so that the geodesic through $(x,\xi)$ is a smoothly closed loop. 

If one wants to go beyond $o(1)$ improvements of~\eqref{e:Linf}, very few results are available. In~\cite{Berard77} (combined with~\cite{Bo16}), the author shows using the Hadamard parametrix that if $(M,g)$ is a manifold without conjugate points, then 
\begin{equation}
\label{e:LinfLog}
\|u_h\|_{L^\infty(M)}\leq C\frac{h^{\frac{1-n}{2}}}{\sqrt{\log h^{-1}}}\|u_h\|_{L^2(M)}.
\end{equation}
The only polynomial improvements that the author is aware of appear in~\cite{I-s} where the authors study Hecke--Maas forms on certain arithmetic surfaces.

\subsection{A review of previous results on averages}
The study of when~\eqref{e:average} is saturated is much more recent and, until the methods of this note were introduced, the only improvements on~\eqref{e:average} available under no additional assumptions on $(M,g)$ are:
\begin{theorem}[\cite{Wym3}]
\label{t:loopAverage}
Suppose that $\sig(\mc{L}_H)=0$. Then
$$
\Big|\int_H u_h d\sigma_H\Big|= o\big(h^{\frac{1-k}{2}}\|u_h\|_{L^2(M)}\big)\
$$
\end{theorem}
The article~\cite{CS} provides $o(1)$ improvements on~\eqref{e:average} on surfaces of negative curvature when $H$ is a geodesic. 

On the other hand, on manifolds with non-positive curvature the Hadamard parametrix is available and as a result logarithmic improvements of the form 
\begin{equation}
\label{e:averageLog}
\Big|\int_Hu_hd\sigma_H\Big|\leq C\frac{h^{\frac{1-k}{2}}}{\sqrt{\log h^{-1}}}\|u_h\|_{L^2(M)}
\end{equation}
hold under a variety of assumptions on the pair $(M,H)$~\cite{SXZ,Wym2,Wym18}. However, none of these results give general dynamical conditions guaranteeing such improvements.

\subsection{Results of the microlocal techniques}

The question raised in all previous attempts to understand when~\eqref{e:Linf} and~\eqref{e:average} can be saturated can be thought of as `In which geometries can saturation occur?'. The question raised in~\cite{GT, CGT,Gdefect,CG17,CG18} is instead `How does an eigenfunction that saturates~\eqref{e:Linf} or~\eqref{e:average} behave?' It then turns out that a sufficiently good understanding of the answer to the latter question yields answers to the former. In fact, by describing the behavior of eigenfunctions saturating~\eqref{e:Linf}, we will be able to extend \emph{all} existing results. Moreover, our analysis of the eigenfunctions saturating these bounds demonstrates that the phenomena governing averages is identical to that governing $L^\infty$ bounds.

\subsubsection{Defect measures}

We begin by describing some of the results of~\cite{CG17} where we rely on defect measures to describe the behavior of $u_h$. Recall that a defect measure is a positive Radon measure, $\mu$, on $\TM$ associated to a sequence of functions $\{u_h\}_{0<h<h_0}$ so that for any $a\in C_c^\infty(T^*M)$, 
$$
\langle Op_h(a)u_h,u_h\rangle_{L^2(M)}\to \int a(x,\xi)d\mu
$$
where $Op_h(a)$ denotes the quantization of the symbol $a$ (see e.g. \cite[Appendix E]{ZwScat} for a description of quantization procedures).
See, for example,~\cite[Chapter 5]{EZB} for a treatment of these measures. We recall that every $L^2$ bounded sequence of functions $\{u_h\}$ has a subsequence with a defect measure, $\mu$ and, moreover if $u_h$ solves~\eqref{e:laplace}, then $\supp \mu\subset \SM$ and $\mu$ is $G^t$ invariant.

 Suppose that $\mu$ is a finite radon measure invariant under the geodesic flow. Then we define for any Borel $A\subset \SNH$,
$$
\mu_H(A):=\lim_{\delta\to 0}\frac{1}{2\delta}\mu\Big(\bigcup_{|t|\leq \delta}G^t(A)\Big).
$$
We write $\omega\perp \nu$ when the measures $\omega$ and $\nu$ are mutually singular.

We then have the following consequence of~\cite[Theorem 6]{CG17} (see also~\cite[Theorem 2]{Gdefect} for the case $k=n$). 
\begin{theorem}
\label{t:micro1}
Let $(M,g)$ be a smooth compact Riemannian manifold and $H\subset M$ a closed embedded submanifold of codimension $k$. Suppose that $u_h$ solves~\eqref{e:laplace} and has defect measure $\mu$. Let $f\in L^1(\SNH;\sig)$ so that 
$$
\mu_H=fd\sig+\omega,\qquad \omega \perp \sig.
$$
Then there is $C_{n,k}>0$ depending only on $(n,k)$ so that for $A\subset H$ with smooth boundary,
$$
\Big|\int_A u_hd\sigma_H\Big|\leq C_{n,k}h^{\frac{1-k}{2}}\int_{\pi_H^{-1}(A)} \sqrt{f}d\sig +o(h^{\frac{1-k}{2}})
$$
where $\pi_H:\SNH\to H$ is the natural projection.
\end{theorem}
Note that Theorem~\ref{t:micro1} can be interpreted as saying that \emph{every} eigenfunction which maximizes either~\eqref{e:Linf} or~\eqref{e:average} must have a component which behaves $o(1)$ microlocally the same as the canonical example on $S^n$. In particular, in order that $u_h$ maximize the $L^\infty$ bounds, there must be a point where $u_h$ behaves like the zonal harmonic, $Z_h$ (See e.g.~\cite[Section 4]{GT} for a description of the defect measures of zonal harmonics.)

As an easy consequence of Theorem~\ref{t:micro1} together with the Poincar\'e recurrence theorem we are able to replace $\mc{L}_H$ in Theorem~\ref{t:loopAverage} by $\mc{R}_H$. 
\begin{corollary}[{{\cite[Theorem 2]{CG17}}}]
\label{c:norecur}
Suppose that $A\subset H$ has smooth boundary, $\sig(\pi_H^{-1}(A)\cap \mc{R}_H)=0$, and $u_h$ solves~\eqref{e:laplace} then
$$
\Big|\int_A u_hd\sigma_H\Big|=o\big(h^{\frac{1-k}{2}}\|u_h\|_{L^2(M)}\big).
$$
\end{corollary}
Using geometric arguments to show that $\sig(\mc{R}_H)=0$ in a variety of settings, we are then able to recover all existing $o(1)$ improvements over~\eqref{e:average} in~\cite[Theorem 4]{CG17}.
\begin{theorem}[{{\cite[Theorem 4]{CG17},\cite[Theorem 3]{CG18}}}]\label{t:app1}
In all of the following situations, $\sig(\mc{R}_H)=0$. 
\begin{enumerate}[label=\textbf{\Alph*.},ref=\Alph*]
\item  \label{a1} $(M,g)$ has no conjugate points and $H$ has codimension $k>\frac{n +1}{2}$. \smallskip 
\item  \label{a2}$(M,g)$ has no conjugate points and $H$ is a geodesic sphere.\smallskip 
\item  \label{a4}$(M,g)$ is a surface with Anosov geodesic flow and $H$ is any submanifold. \smallskip 
\item  \label{a3} $(M,g)$ has  constant negative curvature and $H$ is any submanifold.\smallskip 
\item  \label{a6}$(M,g)$ has Anosov flow, non-positive curvature, and $k>1$.  \smallskip 
\item  \label{a5}$(M,g)$ has Anosov geodesic flow and {non-positive curvature}, and $H$ is totally geodesic.  \smallskip 
\item  \label{a7} $(M,g)$ has {Anosov geodesic flow} and $H$ is a subset $M$ that lifts to a  horosphere. 
\end{enumerate}
\end{theorem}
\smallskip

\subsubsection{Towards quantitative estimates}

In order to pass to the quantitative estimates from~\cite{CG18}, we will first describe some easy consequences of Theorem~\ref{t:micro1}. We say that $A\subset \TM$ is \emph{$[t,T]$ non-self looping} if either 
\begin{equation}
\label{e:nsl}
\begin{gathered}
G^s(A)\cap A=\emptyset,\qquad  s\in[t,T],\qquad
 \text{or}\qquad
G^{-s}(A)\cap A=\emptyset,\qquad  s\in[t,T]\\
\end{gathered}
\end{equation}

We have the following Corollary of Theorem~\ref{t:micro1}.
\begin{corollary}
\label{c:nonLoop}
Suppose that there is a an $h$-independent covering $\big\{B,\{G_\ell\}_\ell\big\}$ of $\SNH$ and $\{t_\ell\}_{\ell}$, $\{T_\ell\}_\ell\subset (1,\infty)$ with $t_\ell<T_\ell$ independent of $h$ so that 
$$
\SNH= B\,\cup \,\bigcup_{\ell} G_\ell.
$$
and $G_\ell$ is $[t_\ell,T_\ell]$ non-self looping an that $u_h$ solves~\eqref{e:laplace}. Then, there is $C>0$ so that for all $u_h$ solving~\eqref{e:laplace}
\begin{equation}
\label{e:prelimEst}
\Big|\int _Hu_hd\sigma_H\Big|\leq Ch^{\frac{1-k}{2}}\left(\sig(B)^{\frac{1}{2}}+\sum_\ell \frac{\sig(G_\ell)^{\frac{1}{2}}t_\ell^\frac{1}{2}}{T_\ell^{\frac{1}{2}}}+o(1)\right)\|u_h\|_{L^2(M)}
\end{equation}
\end{corollary}

In fact, Corollary~\ref{c:norecur} can be deduced from Corollary~\ref{c:nonLoop}. To see this, let $\{U_i\}_{i=1}^\infty$ be a basis for the topology of $\SNH$. Then let $T>0$ and set
$$
E^{\pm,T}_i:=\{x\in U_i\mid G^t(x)\notin U_i, \pm t>T\},\qquad E_i^T=\bigcup_{\pm}E^{\pm,T}_i,\qquad E^\infty_i=\bigcup_{T>0} E_i^T.
$$
Let 
$$
B_{N}=\Big[\bigcap_i \big(\SNH\setminus E^T_i\big)\Big] \cup \Big(\bigcup_{j=N}^\infty E^T_j\setminus \bigcup_{k=1}^{N-1}E_k^T\Big),\qquad \qquad G_i=E^T_i.
$$
Then, since $G_i$ is $[T,S]$ non-self looping for any $S>T$, we apply~\eqref{e:prelimEst} to obtain
$$
\limsup_{h\to 0}\frac{h^{\frac{k-1}{2}}}{\|u_h\|_{L^2(M)}}\Big|\int _Hu_hd\sigma_H\Big|\leq C\Big(\sig(B_N)^{\frac{1}{2}}+\sum_{i=1}^{N-1}\frac{\sig(E^T_i)^{1/2}T^{\frac{1}{2}}}{S^{\frac{1}{2}}}\Big)
$$
Sending $S\to \infty$ gives
$$
\limsup_{h\to 0}\frac{h^{\frac{k-1}{2}}}{\|u_h\|_{L^2(M)}}\Big|\int _Hu_hd\sigma_H\Big|\leq C\sig(B_N)^{\frac{1}{2}}
$$
Sending $N\to \infty$ then gives
$$
\limsup_{h\to 0}\frac{h^{\frac{k-1}{2}}}{\|u_h\|_{L^2(M)}}\Big|\int _Hu_hd\sigma_H\Big|\leq C\sig\Big(\bigcap_i \big(\SNH\setminus E^T_i\big)\Big)^{\frac{1}{2}}
$$
Finally, sending $T\to \infty$ gives
$$
\limsup_{h\to 0}\frac{h^{\frac{k-1}{2}}}{\|u_h\|_{L^2(M)}}\Big|\int _Hu_hd\sigma_H\Big|\leq C\sig\Big(\bigcap_i \big(\SNH\setminus E^\infty_i\big)\Big)^{\frac{1}{2}}
$$

Now, suppose $x$ is not recurrent. Then, there exists $i, T$ so that $x\in U_i$, $G^t(x)\notin U_i$ either for $t>T$ or $-t>T$. In particular, $x\in E_i^T\subset E_i^\infty$. Therefore, if $x$ is not recurrent, then $x\in \cup_iE_i^\infty$. In particular, 
$$
\SNH\setminus \mc{R}_H\,\subset\, \bigcup_i E_i^\infty,\qquad \text{ so }\qquad \bigcap_i \big(\SNH\setminus E_i^\infty\big)\subset \mc{R}_H.
$$
Therefore, if $\mc{R}_H$ has measure 0, then $\bigcap_i\big(\SNH \setminus E_i^\infty\big) $ has measure 0 and we have obtained Corollary~\ref{c:norecur}.

The fact that $\rho\in \mc{R}_H$ does not contain \emph{any} quantitative information about how long it takes for the geodesic through $\rho$ to return to a given neighborhood of $\rho$. Because of this, one should not expect to have a quantitative version of Corollary~\ref{c:norecur}. However, Corollary~\ref{c:nonLoop} \emph{is} quantitative and one might hope that it holds even with $B$, $G_\ell$, and $[t_\ell,T_\ell]$ $h$-dependent. This is almost true, although we will require some additional structure of the sets $B$ and $G_\ell$ (see Theorem~\ref{t:nonLoop2}).

\subsubsection{Quantitative Estimates}

In order to state our quantitative estimates, we will need to define a few additional objects. We will use the metric induced by the Sasaki metric on $\TM$ (see e.g. \cite{Eberlein73} for a description of the Sasaki metric) for convenience, but our results do not depend on the choice of metric on $\TM$. First, fix $\mc{H}\subset \TM$ a smooth hypersurface transverse to the geodesic flow so that $\SNH\subset \mc{H}$.
Define $\psi:\re\times \mc{H}_{\Sigma}\to \TM$ by $\psi(t,q)=\varphi_t(q)$. Next, let
$$
\Tinj:=\sup\{\tau\leq 1: \psi|_{(-\tau,\tau)\times\mc{H}_{\Sigma}}\text{ is injective}\}.
$$

Given $A\subset \TM$, define 
$$
\Lambda_A^\tau:=\bigcup_{|t|\leq \tau}G^t(A).
$$
Then, for $r>0$ and $A\subset \mc{H}$, define 
$$
\Lambda_A^\tau(r):=\Lambda_{A_r}^{\tau+r},\qquad A_r:=\{ \rho\in \mc{H}\mid d(\rho,A)<r\}.
$$
Finally, let $K_H>0$ be a bound for the the sectional curvatures of $H$ and for the second fundamental form of $H$.

\begin{theorem}
\label{t:micro2} 
Let $H\subset M$ be a closed embedded submanifold of codimension $k$. There exist $C_{n,k}>0$ depending only on $n,k$, $\tau_0>0$ depending on $(M,g,\Tinj)$, and $R_0=R_0(n,k,K_H)$ so that the following holds. 

Let $0<\tau<\tau_0$,  $0\leq \class<\frac{1}{2}$, $N>0$,  and $R_0>R(h)\geq 5h^\delta$.
Then, there exists a family $\{\gamma_j\}_{j=1}^{N_h}$ of geodesics through $\SNH$, and a partition of unity $\{\chi_j\}_{j=1}^{N_h}$ for $\LambdaH(h^\class)$ with
 $\chi_j\in S_\class\cap  C^\infty_c(T^*M;[0,1])$,  
 \[
 \supp \chi_j\subset \Tj, \qquad  \qquad \rho_j:=\gamma_j \cap \SNH,
 \]
so that for all $w\in C_c^\infty(H)$, $N>0$ there is $C_N>0$ and $h_0>0$  with the property that for any $0<h<h_0$ and all $u_h$ solving~\eqref{e:laplace}
\begin{align*}
h^{\frac{k-1}{2}}\Big|\int_Hwu_hd\sigma_H\Big|
&\leq C_{n,k} {R(h)^{\frac{n-1}{2}}\sum_{j}\frac{\|Op_h(\chi_j)u_h\|_{L^2(M)}}{\tau^{\frac{1}{2}}}}+C_Nh^N\|u_h\|_{L^2(M)}.
\end{align*}
\end{theorem}
Theorem~\ref{t:micro2} is a much finer analog of Theorem~\ref{t:micro1} and in particular can be interpreted as saying that \emph{every} eigenfunction which maximizes either~\eqref{e:Linf} or~\eqref{e:average} must have a component which behaves the same as the canonical example on $S^n$ microlocally on $h^\delta$ scales. In particular, in order that $u_h$ maximize the $L^\infty$ bounds, there must be a point where $u_h$ behaves like a zonal harmonic at scale $h^\delta$. 

While at first it may seem difficult to use Theorem~\ref{t:micro2} in concrete situations, combining Theorem~\ref{t:micro2} with Egorov's theorem up to the Ehrenfest time (see e.g. \cite{DyGu14}) we obtain a purely dynamical estimate which is readily applicable.

We define the \emph{maximal expansion rate }
$$
\Lambda_{\max}:=\limsup_{|t|\to \infty}\frac{1}{|t|}{\log} \sup_{\SM}\|dG^t(x,\xi)\|.
$$
Then the Ehrenfest time at frequency $h$ is 
$$
T_e(h):=\frac{\log h^{-1}}{2\Lambda_{\max}}.
$$
Note that $\Lambda_{\max}\in[0,\infty)$ and if $\Lambda_{\max}=0$, we may replace it by an arbitrarily small positive constant.  We have the following quantitative version of Corollary~\ref{c:nonLoop}.


\begin{theorem}[{{\cite[Theorem 5]{CG18}}}]
\label{t:nonLoop2}
Suppose that $H\subset M$ is a closed embedded submanifold of codimension $k$. Let $0<\delta<\frac{1}{2}$, $N>0$. There exist positive constants $h_0=h_0(M,g,K_H)$, $\tau_0=\tau_0(M,g,\Tinj)$, $R_0=R_0(n,k,K_H)$ and $C_{n,k}$ depending only on $n$ and $k$, and for each $0<\tau<\tau_0$  there exists and $C_N=C_N(\tau,\delta,M,g)>0$,   so that the following holds.

Let $R_0>R(h)\geq 5h^\delta$,{ $\alpha< 1-2{\limsup_{h\to 0}\frac{\log R(h)}{\log h}}$,} and suppose $\{\Lambda_{_{\rho_j}}^\tau(R(h))\}_{j=1}^{N_h}$  is a cover of $\LambdaH(h^\delta)$ that is the union of $C_{n,k}$ subsets of disjoint tubes (the existence is guaranteed by \cite[Lemma 2.2]{CG18}). In addition, suppose there exist $\mc{B}\subset \{1,\dots N_h\}$ and a finite collection  $\{\mc{G}_\ell\}_{\ell \in L} \subset \{1,\dots N_h\}$ with 
$$
\{1,\dots N(R(h))\}\;\subset\;  \mc{B} \cup \bigcup_{\ell \in L}\mc{G}_\ell,
$$
 and so that for every $\ell \in L$ there exist  $t_\ell(h)>0$ and ${T_\ell(h)}\leq 2 \alpha T_e(h)$  so that
$$
\bigcup_{j\in \mc{G}_\ell}\Lambda_{_{\rho_j}}^\tau(R(h))\;\;\text{ is }\;[t_\ell(h),T_{\ell}(h)]\text{ non-self looping}.
$$
 
Then, for all $w\in C_c^\infty(H)$, $N>0$ there exists $C_N>0$, $h_0>0$ so that for all for $u_h$ solving~\eqref{e:laplace} and $0<h<h_0$,
\begin{multline}
\label{e:estTube}
h^{\frac{k-1}{2}}\Big|\int_H w u_hd\sigma_H\Big|\leq \frac{C_{n,k}\|w\|_{L^\infty}R(h)^{\frac{n-1}{2}}}{\tau^{\frac{1}{2}}}
\!\left[|\mc{B}|^{\frac{1}{2}}+\sum_{\ell \in L }\frac{|\mc{G}_\ell|^{\frac{1}{2}}t_\ell^{\frac{1}{2}}(h)}{T^{\frac{1}{2}}_\ell(h)}+C_Nh^N\right]\!\!\|u_h\|_{L^2(M)}.
\end{multline}
\end{theorem}

Note that the term 
$$
R(h)^{n-1}|\mc{G}_\ell|\propto \sig\Big(\bigcup_{j\in \mc{G}_\ell}\Lambda_{_{\rho_j}}^\tau(R(h))\cap \SNH\Big).
$$
and in particular when $T_{\ell},\,R(h)$ are $h$ independent~\eqref{e:estTube} implies~\eqref{e:prelimEst}. Since~\eqref{e:prelimEst} implies Corollary~\ref{c:norecur}, Theorem~\ref{t:nonLoop2} should be thought of as a quantitative version of the non-recurrent condition. With this intuition in mind, we are able to construct effective covers by tubes in many geometric situations. 
\begin{theorem}[{{\cite[Theorem 3]{CG18}}}]\label{T:applications}
 Let $(M,g)$ be a smooth, compact Riemannian manifold of dimension $n$. Let $H\subset M$ be a closed embedded submanifold of codimension $k$. Suppose one of the following assumptions holds{:}
\begin{enumerate}[label=\textbf{\Alph*.},ref=\Alph*]
\item  \label{a1} $(M,g)$ has no conjugate points and $H$ has codimension $k>\frac{n +1}{2}$. \smallskip 
\item  \label{a2}$(M,g)$ has no conjugate points and $H$ is a geodesic sphere.\smallskip 
\item  \label{a3}$(M,g)$ is a surface with Anosov geodesic flow.  \smallskip 
\item  \label{a6}{$(M,g)$ has Anosov geodesic flow, non-positive curvature and $k>1$.}  \smallskip 
\item  \label{a4}$(M,g)$ has Anosov geodesic flow and {non-positive curvature}, and $H$ is totally geodesic.  \smallskip 
\item  \label{a5} $(M,g)$ has {Anosov geodesic flow} and $H$ is a subset of $M$ that lifts to a horosphere in the universal cover. 
\end{enumerate}
Then there exists $C>0$ so that for all $w\in C_c^\infty(H)$ there is $h_0>0$ so that for $0<h<h_0$ and $u_h$ solving~\eqref{e:laplace}
\begin{equation}
\label{e:subEst}
\Big|\int_Hwu_hd\sigma_H\Big|\leq Ch^{\frac{1-k}{2}}\frac{\|u_h\|_{\LM}}{\sqrt{\log h^{-1}}}.
\end{equation}
\end{theorem}

Finally, there is some uniformity in the estimates from Theorem~\ref{t:micro2} and we can obtain $L^\infty$ estimates. To state these estimates we need to recall a few notions from Riemannian geometry. A Jacobi field along a geodesic $\gamma(t)$ is a vector field along $\gamma(t)$ satisfying
$$
D_t^2J+R(J,\dot{\gamma})\dot \gamma=0
$$
where $D_t$ denotes the covariant derivative along $\gamma$ and $R(\cdot,\cdot)(\cdot)$ denotes the Riemann curvature tensor (see e.g. \cite[Chapter 10]{LeeBook}). We say that $J$ is perpendicular of $\langle J,\dot\gamma\rangle_g=0$ and $\langle D_tJ,\dot\gamma\rangle_g=0$.

For a geodesic $\gamma$, we say that $\gamma$ has a conjugate point of multiplicity $m$ at $t_0$ if there are perpendicular Jacobi fields $\{J_i\}_{i=1}^m $ so that $J_i(0)=0$, $\{D_tJ_i(0)\}_{i=1}^m$ are linearly independent, and $J_i(t_0)=0$. Note that the maximum multiplicity of of a conjugate point is $n-1$ where $n$ is the dimension of the manifold $M$. Moreover, it is not hard to see that there exists $\delta>0$ so that for any geodesic $\gamma$ and any $t_0\in \R$, there are at most $n-1$ conjugate points counted with multiplicity in $[t_0-\delta,t_0+\delta]$.

Define
$$
{\Lambda}_x^{m,r,t}:=\big\{\gamma\in \Lambda: \gamma(0)=x,\,\exists\text{ at least }m\text{ conjugate points to } x \text{ in }\gamma([t-r,t+r])\big\},
$$
where we count conjugate points with multiplicity. Next, for  a set $V \subset M$ write
$$
\mc{C}_{_{\!V}}^{m,r,t}:=\bigcup_{x\in V}\{\gamma(t): \gamma\in \Lambda_x^{m,r,t}\}.
$$
Note that the set $\mc{C}_{_{\!x}}^{n-1,0,t}$ is the set of points that are \emph{maximally conjugate to $x$} at time $t$. In particular, for $y\in \mc{C}_{_{\!x }}^{n-1,0,t}$ there is a geodesic $\gamma$ with $\gamma(0)=x$, $\gamma(t)=y$ and so that \emph{all} of the perpendicular Jacobi fields vanish at $t$. One case where this happens is on the sphere where $x$ and $y$ are antipodal points. While the condition 
$x\notin \mc{C}_x^{n-1,0,t}$ for $t\geq t_0$ is enough to guarantee $o(1)$ improvements in $L^\infty$ bounds, a notion of \emph{uniform maximal self conjugacy} is necessary to have quantitative improvements. 

\begin{theorem}[{{\cite[Theorem 1]{CG18}}}]
\label{t:noConj}
Let $U\subset M$ and suppose that there is $T>0, a>0$ so that for all $x\in U$, 
$$
d\Big(x, \mc{C}_x^{n-1,r_a(t),t}\Big)\geq r_a(t),\qquad t\geq T
$$
where $r_a(t)=a^{-1}e^{-a t}$. Then
$$
\|u_h\|_{L^\infty(U)}\leq C\frac{h^{\frac{1-n}{2}}}{\sqrt{\log h^{-1}}}\|u_h\|_{L^2(M)}.
$$
\end{theorem}
It is not hard to see that Theorem~\ref{t:nonLoop2} implies even stronger estimates where we only assume certain volume control on the directions along which $x$ is maximally self-conjugate. 

We note at this point that Theorems~\ref{T:applications} and~\ref{t:noConj} subsume \emph{all} previous conditions known to give logarithmic improvements and~\ref{t:app1} subsumes \emph{all} previous conditions known to give $o(1)$ improvements.

\bigskip 

\noindent {\sc Acknowledgements.} 
Thanks to Yaiza Canzani for comments on an early version of this note. The author is grateful to the National Science Foundation for support under the Mathematical Sciences Postdoctoral Research Fellowship  DMS-1502661.

\section{The overall ideas of the proofs}

\subsection{The microlocal estimate}

The first important observation in the proof of Theorem~\ref{t:micro2} is that the most localized that an eigenfunction can be is to an $h^{\frac{1}{2}}$ tube around a single length $\sim1$ piece of geodesic.  This is the case, for example, for the highest weight spherical harmonics on $S^2$ given by the restriction of $j^{\frac{n-1}{4}}(x_1+ix_2)^j$ to the sphere, $x_1^2+x_2^2+x_3^2=1$. It is then natural to think of building an eigenfunction out of pieces localized to such tubes. Locally, these pieces are of the form 
$$
u_h(x)=h^{-\frac{1-n}{4}} e^{\frac{i}{h}x_1}e^{-\frac{|x|^2}{2h}}a(x)
$$
with the geodesic given by $\gamma=\{(x_1,0,1,0)\mid |x_1|<1\}$. Here, we have taken $(x_1,x',\xi_1,\xi')$ as coordinates on $\TM$. We will refer to this type of object as a \emph{gaussian beam}.

\begin{figure}
\begin{tikzpicture}
\begin{scope}[shift= {(0,2.3)}]
        \draw [ultra thick] (-4,-1)--(-4,1.5);
        \draw[thick, domain =-5.1:-2.9, samples =100 ]plot (\x, {2.4*e^( -50*(\x+4)*(\x+4))-1}) node[right]{};
                \draw[->](-4,-1.2)-- (-4.3,-1.2);
        \draw[->] (-4,-1.2)--(-3.7,-1.2);
        \node at (-4,-1.5){\scriptsize{$h^{\frac{1}{2}}$}};
        \draw[->] (-5,.5)--(-5,1.4);
        \draw[->] (-5,.5)--(-5,-1);
        \node at (-5.3,.2){\scriptsize{$h^{-\frac{1}{4}}$}};
        \draw [ultra thick,red] (-5.5,-1)--(-2.5,-1);
        \draw node[circle,fill,color=blue,scale=.5] at (-4,-1){};
        \node at (-4,1.9){\scriptsize{Profile across a gaussian beam}};
        \node at (-4,-2.3){\scriptsize{Profile along a gaussian beam}};
\end{scope}

    \draw[ domain=-5.5:-2.5,samples=500]   plot (\x,{sin( (15 *\x)  r )-1.5})   node[right] {};
    
     \draw[->](-4,-2.7)-- (-2.5,-2.7);
        \draw[->] (-4,-2.7)--(-5.5,-2.7);
        \node at (-4,-2.9){\scriptsize{$Ch$}};
        \draw[->] (-6,-1.5)--(-6,-.5);
        \draw[->] (-6,-1.5)--(-6,-2.5);
        \node at (-6.5,-1.3){\scriptsize{$h^{-\frac{1}{4}}$}};
        \draw[->,blue] (-5.5,-.3)--(-2.5,-.3);
        
        \begin{scope}[shift={(0,2)}]

        \draw[ultra thick, red] (0,1.5)-- (0,-1.5);
        \node at (0,1.7){\footnotesize{{$H$}}};
        \foreach \t in {0}
        { \draw[->,blue] (-1.5, \t) -- (1.5,\t); }
\node at (3.3,1.2){$\|u_h\|_{L^2}=1$} ;
\node at (3.3,.2){$\int_{H}u_hd\sigma_{H}\sim ch^{\frac{1}{4}}$} ;
\end{scope}

\begin{scope}[shift={(6,0)}]
\draw[ domain=-5.5:-2.5,samples=500]   plot (\x,{e^( -10*(\x+4)*(\x+4))-1.5})   node[right] {};
    
     \draw[->](-4,-2.7)-- (-2.5,-2.7);
        \draw[->] (-4,-2.7)--(-5.5,-2.7);
        \node at (-4,-2.9){\scriptsize{$Ch^{\frac{1}{2}}$}};
        \draw[->] (-6,-1.5)--(-6,-.5);
        \draw[->] (-6,-1.5)--(-6,-2.5);
        \node at (-6.5,-1.3){\scriptsize{$h^{-\frac{1}{4}}$}};
        \draw[->,white] (-5.5,-.3)--(-2.5,-.3)node[above,midway,black]{\scriptsize{Profile after restriction to $H$}};
\end{scope}
\end{tikzpicture}
\caption{\label{f:normal} Diagram when a gaussian beam passes over $H$ normally}
\end{figure}
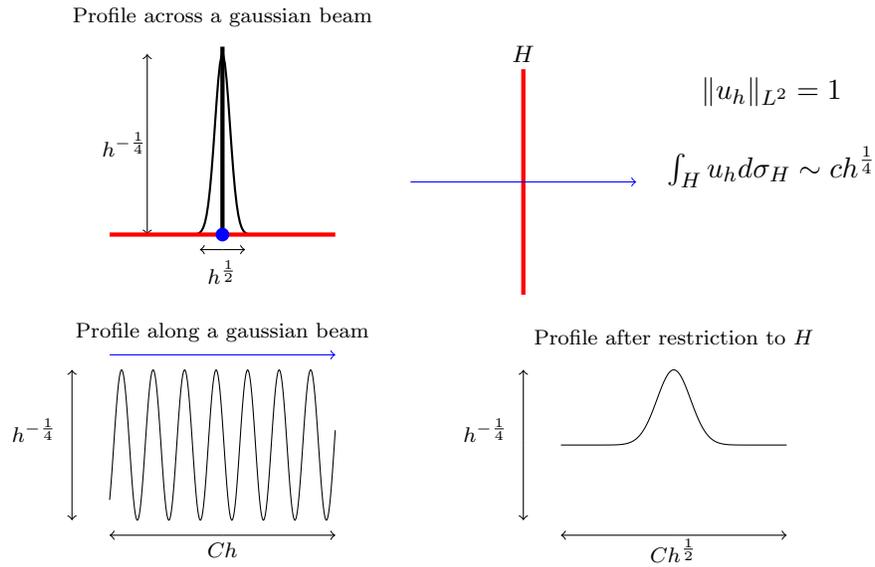

The first step is then to understand how the average over $H$ of an eigenfunction localized to such a tube behaves. This is a two step process. First, if the tube is passing over the hypersurface in a direction which is not normal to the hypersurface, then the contribution is $O(h^\infty)$. Such a restriction is shown in Figure~\ref{f:nonnormal}. Since oscillation remains after restriction, the contribution from such a tube is $O(h^\infty)$. Once we have this in place, we need to study tubes passing normally over $H$ as in Figure~\ref{f:normal}.

\begin{figure}
\begin{tikzpicture}
\begin{scope}[shift= {(0,2.3)}]
        \draw [ultra thick] (-4,-1)--(-4,1.5);
        \draw[thick, domain =-5.1:-2.9, samples =100 ]plot (\x, {2.4*e^( -50*(\x+4)*(\x+4))-1}) node[right]{};
                \draw[->](-4,-1.2)-- (-4.3,-1.2);
        \draw[->] (-4,-1.2)--(-3.7,-1.2);
        \node at (-4,-1.5){\scriptsize{$h^{\frac{1}{2}}$}};
        \draw[->] (-5,.5)--(-5,1.4);
        \draw[->] (-5,.5)--(-5,-1);
        \node at (-5.3,.2){\scriptsize{$h^{-\frac{1}{4}}$}};
        \draw [ultra thick,red] (-5.5,-1)--(-2.5,-1);
        \draw node[circle,fill,color=blue,scale=.5] at (-4,-1){};
        \node at (-4,1.9){\scriptsize{Profile across a gaussian beam}};
        \node at (-4,-2.3){\scriptsize{Profile along a gaussian beam}};
\end{scope}

    \draw[ domain=-5.5:-2.5,samples=500]   plot (\x,{sin( (15 *\x)  r )-1.5})   node[right] {};
    
     \draw[->](-4,-2.7)-- (-2.5,-2.7);
        \draw[->] (-4,-2.7)--(-5.5,-2.7);
        \node at (-4,-2.9){\scriptsize{$Ch$}};
        \draw[->] (-6,-1.5)--(-6,-.5);
        \draw[->] (-6,-1.5)--(-6,-2.5);
        \node at (-6.5,-1.3){\scriptsize{$h^{-\frac{1}{4}}$}};
        \draw[->,blue] (-5.5,-.3)--(-2.5,-.3);
        
        \begin{scope}[shift={(0,2)}]

        \draw[ultra thick, red] (0,1.5)-- (0,-1.5);
        \node at (0,1.8){\footnotesize{$H$}};
        \draw[->,blue] (45:-1.5)--(45:1.5);
\node at (3.3,1.2){$\|u_h\|_{L^2}=1$} ;
\node at (3.3,.2){$\int_{H}u_hd\sigma_{H}=O(h^\infty)$} ;
\end{scope}

\begin{scope}[shift={(6,0)}]
\draw[ domain=-5.5:-2.5,samples=500]   plot (\x,{e^( -10*(\x+4)*(\x+4))*sin((40*\x) r)-1.5})   node[right] {};
    
     \draw[->](-4,-2.7)-- (-2.5,-2.7);
        \draw[->] (-4,-2.7)--(-5.5,-2.7);
        \node at (-4,-2.9){\scriptsize{$Ch^{\frac{1}{2}}$}};
        \draw[->] (-6,-1.5)--(-6,-.5);
        \draw[->] (-6,-1.5)--(-6,-2.5);
        \node at (-6.5,-1.3){\scriptsize{$h^{-\frac{1}{4}}$}};
        \draw[->,white] (-5.5,-.3)--(-2.5,-.3)node[above,midway,black]{\scriptsize{Profile after restriction to $H$}};
\end{scope}
\end{tikzpicture}
\caption{\label{f:nonnormal} Diagram when a gaussian beam passes over $H$ non-normally}
\end{figure}
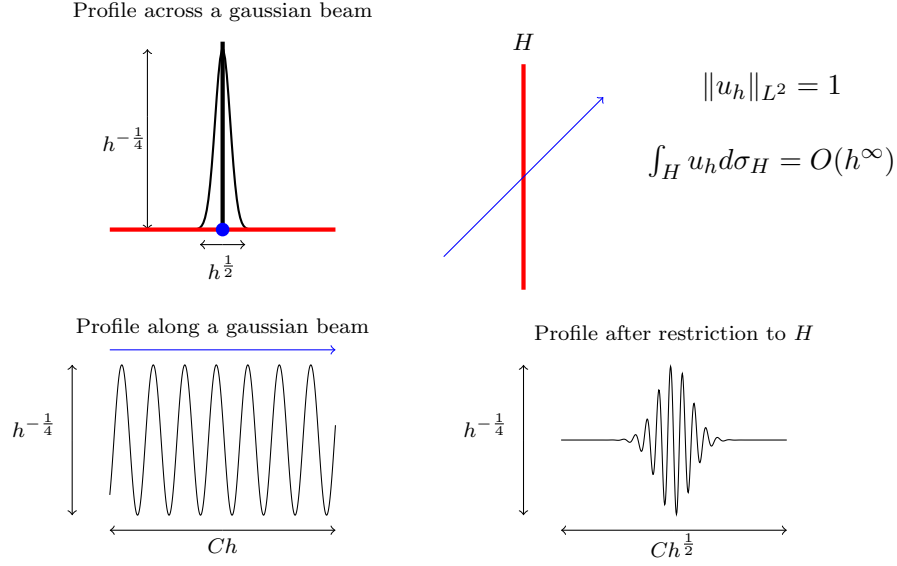

When we decompose eigenfunctions using tubes, we will use tubes of size $R(h)\geq h^\delta$ for some $0\leq\delta<\frac{1}{2}$ so that the symbolic calculus works well. It remains to understand how eigenfunctions localized to such tubes behave when restricted to submanifolds. The key observation is that localization to a small tube implies better control on oscillation and that this control gives improved $L^\infty$ estimates. In particular, imagine that we are working on $\R^n$ with coordinates $(x_1,x',\xi_1,\xi')$ on $T^*\R^n$ and
$$
\gamma=\{(x_1,0,1,0)\mid |x_1|<1\}\subset T^*\R^n.
$$
Then, assume that a function $u_h$ has frequencies only in $|\xi'|\leq R(h)$, i.e. with
$$
\mc{F}_h(u_h)(\xi):=\int e^{-\frac{i}{h}\langle x,\xi\rangle} u_h(x)dx,
$$
satisfying $\supp \mc{F}_h(u)\subset \{|\xi'|\leq R(h)\}$ modulo $O(h^\infty)$. Then
$$
\|(hD_{x'})^mu_h\|_{L^2}\leq CR(h)^m\|u_h\|_{L^2}.
$$
In particular, $u$ is oscillating at frequency $R(h)h^{-1}$ in the $x'$ variables. 

Now, suppose that $H$ is given by $\{x_1=\bar{x}=0\}$ where $\bar{x}\in \R^{k-1}$ and $x'=(\bar{x},x'')$. By the standard Sobolev embedding $H^s(\R^{m})\to L^\infty(\R^{k-1})$ for $s>\frac{k-1}{2}$, such oscillation then implies that 
$$
\|u_h(x_1,\cdot)\|_{L^\infty_{\bar{x}}L^2_{x''}}\leq CR(h)^{\frac{k-1}{2}}h^{\frac{1-k}{2}}\|u_h(x_1,\cdot)\|_{L^2_{x'}}.
$$
Then, if we assume in addition $(-h^2\Delta-1)u=0$, a standard energy estimate (see e.g. \cite[Chapter 7]{EZB}) implies that for $u$ localized close enough  ($h$ independently) to $\gamma$, 
$$
\|u_h(x_1,\cdot)\|_{L^2_{x'}}\leq C\|u_h\|_{L^2}.
$$
In particular,
$$
\|u_h\|_{L^\infty_{\bar{x}}L^2_{x''}}\leq CR(h)^{\frac{k-1}{2}}h^{\frac{1-k}{2}}\|u_h\|_{L^2}.
$$

Finally, if $u$ is also supported on $|x'|\leq R(h)$ modulo $O(h^\infty)$, then
\begin{equation}
\label{e:singleTube}
\Big|\int u_h(0,x'')dx''\Big|\leq CR(h)^{\frac{k}{2}}\|u_h(0,\cdot)\|_{L^2_{x''}}\leq CR(h)^{\frac{n-1}{2}}h^{\frac{1-k}{2}}\|u_h\|_{L^2}.
\end{equation}

\begin{figure}
\begin{tikzpicture}
\draw[thick, dashed, ->, scale=1.3] (-{sqrt(3)},-.5)-- ({sqrt(3)},.5)node[right] {$\gamma$};
\draw[thick] (-{sqrt(3)},-.5) ellipse (.2 and .375);
\draw[thick] ({sqrt(3)},.5) ellipse (.2 and .375);
\draw[thick] (-{sqrt(3)},-.5+.375)-- ({sqrt(3)},.5+.375);
\draw[thick] (-{sqrt(3)},-.5-.375)-- ({sqrt(3)},.5-.375);
\draw[gray, shift={(-.5,-.75)}] (-30:0)-- (-30:1.5)--++(90:2.5)--++(150:1.5) node[above]{$S\!N^*\!H$}--cycle; 
\draw (0,0) node[below]{$\rho$};
\fill (0,0) circle (.05);
\end{tikzpicture}
\caption{\label{f:tube} A single tube}
\end{figure}
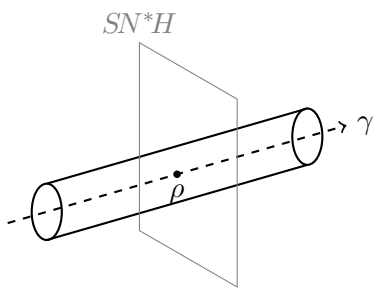

In order to make this argument on a general manifold, we construct microlocal cutoffs, $\chi$, to $R(h)$ sized tubes around geodesics (see Figure~\ref{f:tube}) which essentially commute with the Laplacian near $H$. We are then able to use the calculus of pseudodifferential operators to obtain the estimate~\eqref{e:singleTube}.

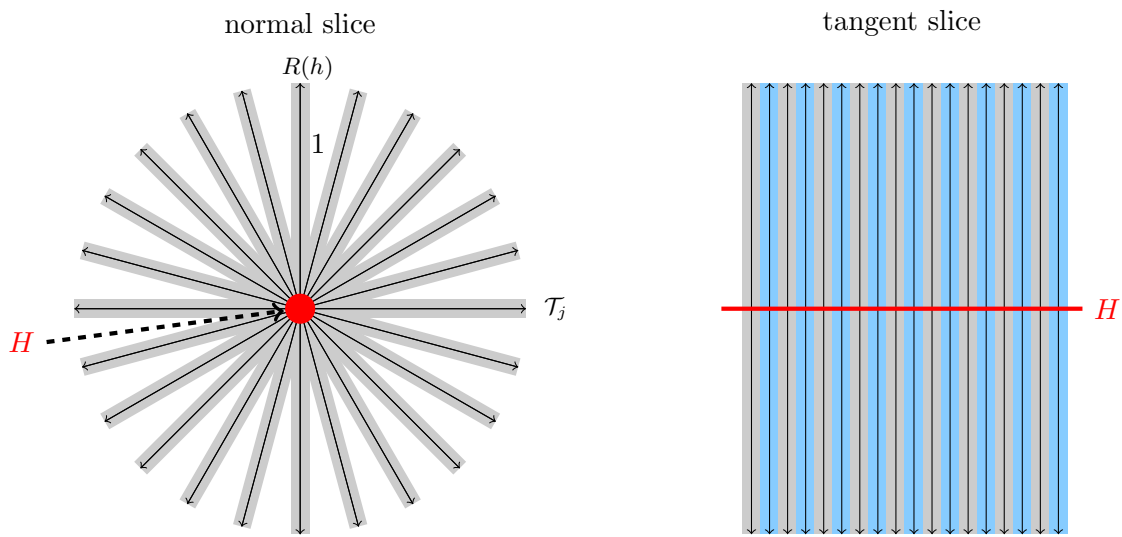
\begin{figure}
\begin{tikzpicture}  
\begin{scope}[scale=2]
        \foreach \t in {0,15,..., 359}
        { 
        \fill[light-gray,rotate=\t] (-.06,0) rectangle (.06,1.5);}    
        \foreach \t in {0,15,..., 359}
        { 
        \draw[->] (\t:-1.5) -- (\t:1.5); }       
\node at (.05,1.6){\footnotesize{$R(h)$}} ;
\node at (.12,1.1){$1$} ;
\node at (1.7,0){\footnotesize{{$\mathcal{T}_j$}}};
\fill[red](0,0)circle (.1); 
         \draw[ultra thick,dashed,->] (7.5:-1.7)node[left]{\red{$H$}}--(7.5:-.11) ;
         \node at (0,1.9){normal slice};
\end{scope}

\begin{scope}[scale=2,shift={(3,0)}]
        \foreach \t in {0,.24,..., 2}
        { 
        \fill[fill=light-gray,shift={(\t,0)}] (-.06,0) rectangle (.06,1.5);
        \fill[fill=light-blue,shift={(\t+.12,0)}](-.06,0) rectangle (.06,1.5);
        \draw[->,shift={(\t,0)}] (0,0) -- (0,1.5);
        \draw[->,shift={(\t+.12,0)}] (0,0) -- (0,1.5);}        
        \foreach \t in {0,.24,..., 2}
        { 
        \fill[light-gray,shift={({\t},0)}] (-.06,0) rectangle (.06,-1.5);
        \fill[fill=light-blue,shift={(\t+.12,0)}](-.06,0) rectangle (.06,-1.5);
        \draw[->,shift={({\t+.12},0)}] (0,0) -- (0,-1.5);
        \draw[->,shift={({\t},0)}] (0,0) -- (0,-1.5);}       
\draw[ultra thick,red] (-.2,0)--(2.2,0) node[right]{$H$};
\node at (1,1.9){tangent slice};
\end{scope}
\end{tikzpicture}
\caption{\label{f:cover} The cover of $\SNH$ by tubes. (left) The projection onto a plane normal to $H$. (right) The projection of the tubes onto a plane tangent to $H$. Note that each pair of tubes (pointing up and down) on the right corresponds to a whole sphere of tubes in $\SNH$ as pictured on the left. The tubes alternate in color only to make it easier to distinguish adjacent tubes.}
\end{figure}

In order to complete the proof of Theorem~\ref{t:micro2}, we then cover $\SNH$ by tubes as in Figure~\ref{f:cover}. In the case of $k=n$, combining the estimates is just a matter of applying the triangle inequality. However, when $k<n$, we must once again use that, modulo $O(h^\infty)$, the cutoffs are supported \emph{in physical space} at a distance $R(h)$ from a geodesic. Covering $H$ by balls of radius $R(h)$ and applying the triangle inequality in each ball then gives the required estimate.

\subsection{From the microlocal estimate to a Theorem~\ref{t:nonLoop2}}

Passing from Theorem~\ref{t:micro2} to Theorem~\ref{t:nonLoop2} is an application of Egorov's theorem to long times. In particular, observe that if $\chi\in C_c^\infty(\TM;[0,1])\cap S_\delta$ is $[t_0,T_0]$ non-self looping, then 
\begin{equation}
\label{e:propSymb}
\Big|\frac{1}{T}\int_0^T\chi^2\circ\varphi_tdt\Big|\leq \frac{t_0}{T}.
\end{equation}
Here, we say $\chi\in S_\delta$ if 
$$
\big|\partial_x^\alpha\partial_{\xi}^\beta \chi(x,\xi)\big|\leq C_{\alpha\beta}h^{-\delta(|\alpha|+|\beta|)}.
$$
Together with Egorov's theorem to the Ehrenfest time (see e.g. \cite[Proposition 3.8]{DyGu14}),~\eqref{e:propSymb} implies that 
\begin{align*}
\|Op_h(\chi)u_h\|_{L^2}^2&=\left\langle \frac{1}{T}\int_0^{T}e^{ith\Delta_g}Op_h(\chi)^*Op_h(\chi)e^{-ith\Delta_g}dtu_h,u_h\right\rangle _{L^2(M)}\\
&\leq \frac{t_0}{T}\big(1+O(h^\e)\big)\|u_h\|_{L^2(M)}^2.
\end{align*}

In particular, if $\bigcup_{j\in \mc{G}_\ell}\Tj$ is $[t_\ell,T_\ell]$ non-self looping, then since there are at most $C$ $\chi_j$ with overlapping supports,
\begin{equation}
\label{e:nonLoopingTubesEst}
\sum_{j\in \mc{G}_\ell}\|Op_h(\chi_j) u_h\|^2_{L^2}\leq C\Big\|\sum_{j\in \mc{G}_\ell}Op_h(\chi_j)u_h\Big\|_{L^2}^2\leq \frac{t_\ell}{T_\ell}\big(1+O(h^\e)\big)\|u_h\|_{L^2(M)}^2.
\end{equation}
An application of Cauchy--Schwarz together with Theorem~\ref{t:micro2} and~\eqref{e:nonLoopingTubesEst} then gives Theorem~\ref{t:nonLoop2}.

\subsection{Construction of effective covers}

There are two mechanisms used to construct the effective covers for Theorem~\ref{T:applications}; contraction and rotation. 

\subsubsection{Contraction} 
\begin{figure}[h]
\begin{tikzpicture}
\begin{scope}[scale=4]
 \draw[color=blue,thick, domain =0:2, samples =100 ]plot (\x, {.5*e^( -1*(\x)}) node[right]{};
 \draw[color=blue,thick, domain =0:2, samples =100 ]plot (\x, -{.5*e^( -1*(\x)}) node[right]{};

\foreach \y in{0,0.5,...,2}{
\draw[red,thick] (\y, -.7)-- (\y, .7)node[above]{\tiny{$H$}} ;
\draw[ultra thick] (\y, -{.5*e^( -1*(\y)})-- (\y, {.5*e^( -1*(\y)}) node[midway,left]{\tiny{$G^t(B)$}};
}
\draw[thick,->] (0,-.75)--(2,-.75)node[right]{\tiny{$t\to\infty$}};
\end{scope}
\end{tikzpicture}
\caption{\label{f:tubeContract} Contraction mechanism for constructing effective covers.}
\end{figure}
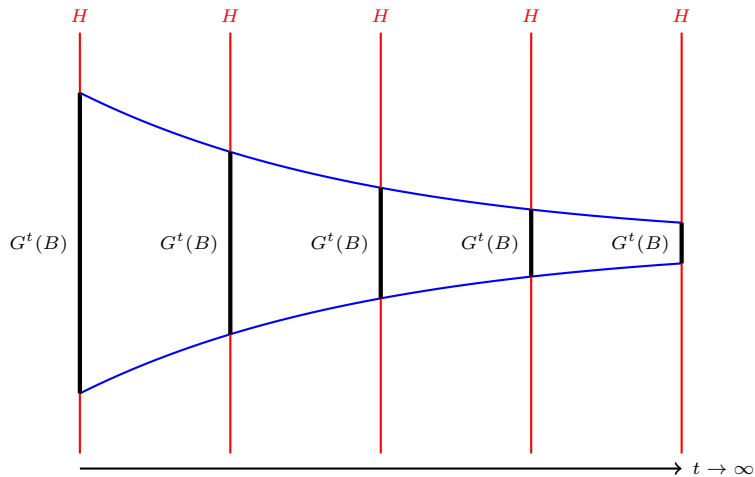

In the contraction mechanism, pictured in Figure~\ref{f:tubeContract}, we use the fact that a subset of $\SNH$ contracts under the flow either forward or backward in time. This is the case, for instance, when $H$ is contained in a stable or unstable horosphere. Under this condition, we start with a macroscopic set $A_0$ and flow it forward in time. We remove all of the pieces of
$$
B_0:=\bigcup_{t_0}^TG^t(A_0)
$$
intersecting $A_0$. Since $A_0$ is contracting, we may choose $t_0$ large enough so that $\sig(B_0)\leq \e\sig(A_0)$. We then let $G_0=A_0\setminus B_0$ and $A_1=B_0$. By construction $G_0$ is $[t_0,T]$ non-self looping. We can then repeat the process replacing $T$ by $\frac{1}{2}T$ to obtain $G_1$ that is $[t_0,\frac{1}{2}T]$ non-self looping. Inductively repeating this process, we construct an effective non-self looping cover of $\SNH$.

\subsubsection{Rotation}
\begin{figure}[h]
\begin{tikzpicture}
\foreach \t in{0}{
\begin{scope}[shift= {(2.4*\t-1.5,-4)},scale=1.5]
\draw[thick,->](\t,-1.7)node[below]{\tiny{${t=\t}$}}--(\t,1.7);
\draw[gray,->](\t,0)--(\t+.5,.5)node[right]{\color{gray}{\tiny{$H_{|\xi|_g^2}$}}};
\draw[thick,dashed,blue](\t-.3,-.9)--(\t+.4,1.2)node[right]{\tiny{\color{blue}$T_\rho \SNH$}};
\draw[thick,->](\t-.85,.5)--(\t+.85,-.5);
\draw[thick,red] ({\t-.75*.2/(\t+1)},{-.75*.6*(\t+1)}) --({\t+.75*.2/(\t+1)}, {.75*.6*(\t+1)})node[left]{\tiny{\color{red}$G^t(B_\rho)$}};
\end{scope}
}
\foreach \t in{1,...,2}{
\begin{scope}[shift= {(2.4*\t-1.5,-4)},scale=1.5]
\draw[thick,->](\t,-1.7)node[below]{\tiny{${t=\t}$}}--(\t,1.7);
\draw[gray,->](\t,0)--(\t+.5,.5)node[right]{\color{gray}{\tiny{$H_{|\xi|_g^2}$}}};
\draw[thick,dashed,blue](\t-.3,-.9)--(\t+.4,1.2)node[right]{\tiny{\color{blue}$T_{_{\!G^{\t}(\rho)}} \SNH$}};
\draw[thick,->](\t-.85,.5)--(\t+.85,-.5);
\draw[thick,red] ({\t-.75*.2/(\t+1)},{-.75*.6*(\t+1)}) --({\t+.75*.2/(\t+1)}, {.75*.6*(\t+1)})node[left]{\tiny{\color{red}$G^t(B_\rho)$}};
\end{scope}
}
\end{tikzpicture}
\caption{ \label{f:tubeRotate} Rotation mechanism for constructing effective covers.}
\end{figure}
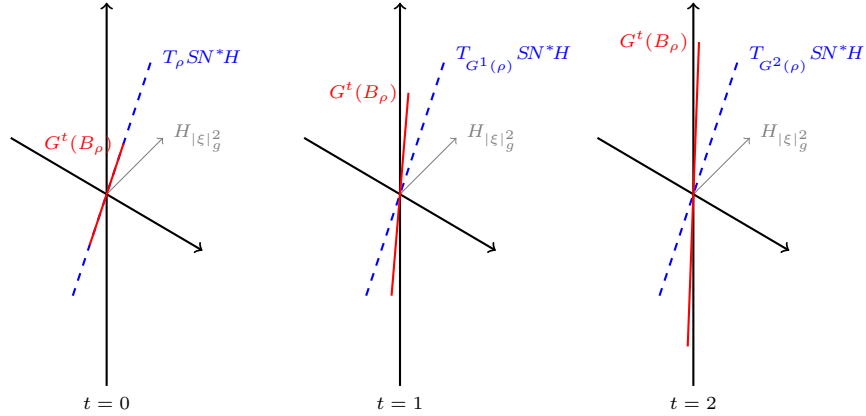

In the rotation mechanism, pictured in Figure~\ref{f:tubeRotate}, a ball of small radius $B_\rho\subset\SNH$ rotates to become transverse to the plane of $T_{G^t(\rho)}\SNH$ (when $G^t(\rho)\in \SNH$) as $t\to \pm \infty$. In this situations, we can use the implicit function theorem to show that the intersection of $\bigcup_{t_0}^T\varphi_t(B_\rho)$ with $\SNH$ is a finite union of lower dimensional subsets. Covering these lower dimensional subsets by tubes with small volume, we are able to construct an effective cover.

\subsubsection{Effective covers with no uniformly maximal self-conjugate points}
\begin{figure}[h]
\includegraphics[width=\textwidth]{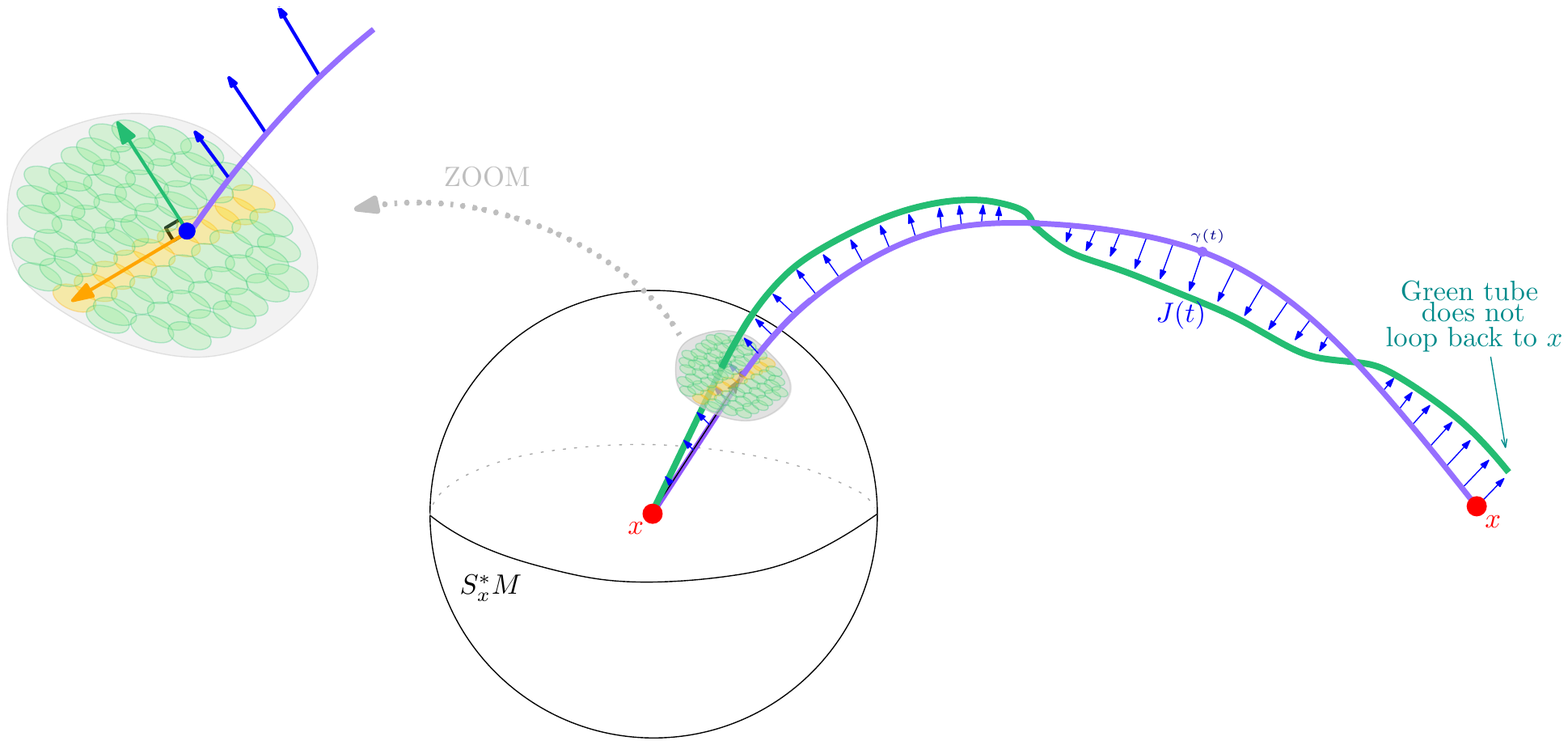}
\caption{\label{f:noConj} The construction of a non-self looping collection of tubes under a non-uniformly maximal self conjugacy assumption.}
\end{figure}

We now sketch the construction of the cover by tubes that is used to prove Theorem~\ref{t:noConj}. The crucial fact is that if $J$ is a Jacobi field along $\gamma$ with $J(0)=0$ and $\Gamma:(-\e,\e)\times \R\to M$ is a map so that 
$$
\begin{gathered}
\Gamma(s,0)=\gamma(0),\qquad\qquad \gamma_s:t\mapsto\Gamma(s,t)\text{ is a geodesic},\qquad\qquad\partial_sD_t\Gamma(0,0)=D_tJ(0),
\end{gathered}
$$
then $\partial_s\gamma(0,t)=J(t)$. Said another way, if $J(t_0)$ is non-zero, then for $s\neq 0$ small $\gamma_s(t_0)\neq \gamma(0)$. 

Translating this from the $SM\subset TM$ to $\SM\subset \TM$, this implies that there is a vector $V=(D_tJ(0))^\sharp\in T_{{\dot\gamma}^\sharp}S^*_xM$ so that
$$
d\pi dG^{t_0} V\neq 0
$$
where $\pi:\TM\to M$ denotes the projection. Using this together with the implicit function theorem, we find a submanifold $B\subset S^*_xM$ of dimension$<n-1$ and a neighborhood $W$ of ${\dot\gamma}^\sharp$ so that for $t$ near $t_0$ and $\rho\in W\setminus B$, $G^t(\rho)\notin S^*_xM$. We can then cover $B$ by $\sim R(h)^{2-n}$ tubes. 

Since $x$ is not maximally self-conjugate for $t>s_0$, we can repeat this argument near each point $\rho\in  S^*_xM$ and then for approximately $T$ values of $t_0$, we produce a large collection of tubes, $\mc{G}$ whose union is $[s_0,T]$ non-self looping and $\sim T R(h)^{2-n}$ possibly looping tubes $\mc{B}$.

In order to make this construction work, we must control the size of the neighborhood $W$ near each $\rho$. It is precisely in this quantification where the uniformity in the non-maximally self conjugacy is used.

\bibliography{biblio}
\bibliographystyle{alpha}

\end{document}